\documentclass{birkjour}
 \newtheorem{thm}{Theorem}[section]

 \newtheorem{prop}[thm]{Proposition}
 \theoremstyle{definition}
 \newtheorem{defn}[thm]{Definition}
 \theoremstyle{remark}
 \newtheorem{rem}[thm]{Remark}
 \newtheorem*{ex}{Example}
 \numberwithin{equation}{section}

\DeclareMathOperator{\Tr}{Tr}
\DeclareMathOperator{\Disc}{Disc}

\newcommand{\anelok}{\mathcal{O}_{\mathbb{K}}} 

\usepackage{xcolor}
\usepackage{comment}

\begin{document}

%
%
%
%
%
%
%
%
%

\title[Constructions of lattices via
biquadratic and triquadratic fields ]
 {Constructions of lattices via
biquadratic and triquadratic fields}

\author[L. F. Santos and G. C. Jorge]{L. F. Santos and G. C. Jorge}

\address{%
Institute of Science and Technology \\
Federal University of São Paulo (UNIFESP) \\
São José dos Campos - Brazil}

\email{leonardo.farias@unifesp.br }
\email{grasiele.jorge@unifesp.br}

\thanks{This work was partially supported by FAPESP 2023/05215-7, CAPES, PICME and CNPq 405842/2023-6.}




\begin{abstract} The lattices $D_4$ and $E_8$ are known to be the densest lattices in dimensions 4 and 8, respectively. In this paper, we employ tools from algebraic number theory to prove that the $D_4$-lattice arises from an infinite family of totally imaginary biquadratic fields. Furthermore, we extend this construction to show that both the $D_8$ and $E_8$ can be realized via triquadratic fields derived from this family.
\end{abstract}

\maketitle

\section{Introduction}	

Lattices are discrete sets of points in the \( n \)-dimensional Euclidean space \( \mathbb{R}^n \), which are described as all integer linear combinations of independent vectors. A sphere packing in \( \mathbb{R}^n \) is a distribution of spheres of equal radius such that any two spheres either do not intersect or intersect only on their boundaries. A lattice packing is a sphere packing in which the set of centers of the spheres forms a lattice. The packing density of a sphere packing is the proportion of \( \mathbb{R}^n \) covered by the union of the spheres of the packing. Usually, the problem of finding good signal constellations for a Gaussian channel is associated with the search for lattices with high packing density \cite{Conway}. 

Given a number field \( \mathbb{K} \) of degree \( n \), one can construct lattices in \( \mathbb{R}^n \) from \( \mathbb{Z} \)-modules contained in \( \mathbb{K} \), using either the canonical or the twisted embedding \cite{Eva2,Eva,Samuel}. A lattice obtained in this way is called an algebraic lattice. An interesting problem, both from a theoretical and a practical perspective, is to characterize which lattices can arise from certain number fields. In this paper, we focus on biquadratic and triquadratic fields, i.e., number fields of the forms $\mathbb{Q}(\sqrt{a},\sqrt{b})$ and $\mathbb {Q}(\sqrt{a},\sqrt{b},\sqrt{c}),$ respectively, where $a$, $b$ and $c$ are square-free integers.

The $D_n$-lattice is defined as the set of points in \( \mathbb{R}^n \) with integer coordinates whose sum is even. The lattice packing associated with \( D_4 \) has been proven to be the densest lattice packing in \( \mathbb{R}^4 \), with density \( \pi^2/16 \) (see \cite{Korkine}). The \( E_8 \)-lattice is defined as the set of points in \( \mathbb{Z}^8 \cup (\mathbb{Z} + \tfrac{1}{2})^8 \) whose coordinate sum is even. In \cite{Viazovska} it was proved that the sphere packing (not just the lattice packing) associated with \( E_8 \) achieves the highest possible density in dimension 8. This density is equal to \( \pi^4 / 384 \), which is precisely twice that of the packing associated with \( D_8 \), whose density is \( \pi^4 / 768 \).

Algebraic constructions of \( D_4 \) have been considered in some papers. The field \( \mathbb{Q}(\zeta_8) = \mathbb{Q}(\sqrt{2}, i) \) is used in \cite{Andrade, Eva, Boutros}, while the fields \( \mathbb{Q}(\zeta_{15} + \zeta_{15}^{-1}) \), \( \mathbb{Q}(\zeta_{16} + \zeta_{16}^{-1}) \)  and \( \mathbb{Q}(\zeta_{12}) = \mathbb{Q}(\sqrt{3}, i) \) appear in \cite{Jorge, Jorge2, Jorge-Densos}. An infinite family of totally real biquadratic fields of the form \( \mathbb{Q}(\sqrt{2}, \sqrt{q}) \), where \( q \equiv 3 \pmod{8} \) is prime, yielding lattices $D_4$, is presented in \cite{Interlando}. On the other hand, the case of \( D_8 \) is treated in \cite{Jorge, Jorge2}, based on the fields \( \mathbb{Q}(\zeta_{17} + \zeta_{17}^{-1}) \) and \( \mathbb{Q}(\zeta_{32} + \zeta_{32}^{-1}) \).

As for \( E_8 \), several constructions have been proposed. The fields \( \mathbb{Q}(\zeta_{20})\) and \( \mathbb{Q}(\zeta_{24}) = \mathbb{Q}(\sqrt{2}, \sqrt{3}, i) \) were used in \cite{Craig}, while \( \mathbb{Q}(\zeta_{15}) \), \( \mathbb{Q}(\zeta_{20}) \) and \( \mathbb{Q}(\zeta_{24}) \) are employed in \cite{Eva}. The field \( \mathbb{Q}(\zeta_{20}) \) is also considered in \cite{Andrade}. In \cite{Trajano}, it is shown that there are infinitely many subfields \( \mathbb{K} \subseteq \mathbb{Q}(\zeta_{pq}) \), with \( p \) and \( q \) distinct primes, that can generate rotated \( E_8 \)-lattices. Further constructions involving the field \( \mathbb{Q}(\zeta_{60} + \zeta_{60}^{-1}) = \mathbb{Q}(\sqrt{2},\sqrt{3}, \sqrt{5})\) are presented in \cite{Eva-Ivan, Jorge-Densos}.

In this paper, we focus on constructions of $D_4$, $D_8$ and $E_8$ via biquadratic and triquadratic fields. In Section~3, we show that the $D_4$-lattices can be obtained via an infinite family of totally imaginary biquadratic fields of the form \( \mathbb{Q}(\sqrt{2}, \sqrt{-q}) \), where $q \equiv 3 \pmod 8$ is a product of primes congruent to either $1$ or $3$ modulo $8$.  In Section~4, we describe methods for constructing the lattices \( D_4 \oplus D_4 \), \( D_8 \), and \( E_8 \) via the triquadratic fields $\mathbb{Q}(\sqrt{2}, \sqrt{q},i)$.

\section{Preliminaries}
This section presents some elementary results on algebraic lattices. More details can be found in \cite{Conway,Samuel, Stewart}. Let $\Lambda \subset \mathbb{R}^n$ be a full-rank lattice and $\rho$ the packing radius of $\Lambda$. The packing density of a lattice is given by
\[
\Delta(\Lambda) = 
\frac{
\text{Volume of the ball of radius } \rho
}{
\text{Volume of a fundamental region of } \Lambda
}
= 
\frac{\mathrm{vol}(B[0,1]) \, \rho^n}{\mathrm{vol}(\Lambda)}.
\]
Since for fixed $n$ the volume of $B[0,1]$ is constant, the center density of a lattice is defined as
\[
\delta(\Lambda) = \frac{\rho^n}{\mathrm{vol}(\Lambda)}.
\]
Let \( \mathbb{K} \) be a CM-field of degree \( n = 2r \), and let \( \sigma_1, \sigma_2, \ldots, \sigma_n \) denote the \( n \) embeddings of \( \mathbb{K} \) into \( \mathbb{C} \), ordered such that \( \overline{\sigma_{j+r}} = \sigma_j \), for all \( j = 1, \ldots, r \), where \( \overline{\sigma_{j+r}} \) denotes the composition of \( \sigma_{j+r} \) with complex conjugation.
 The canonical embedding \( \sigma_{\mathbb{K}} : \mathbb{K} \to \mathbb{R}^n \) is defined, for each \( x \in \mathbb{K} \), by
\begin{align*}
    \sigma_{\mathbb{K}}(x) &= \left(\Re(\sigma_1(x)), \Im(\sigma_1(x)), \ldots, \Re(\sigma_r(x)), \Im(\sigma_r(x)) \right),
\end{align*}
where \( \Re \) and \( \Im \) denote the real and imaginary parts of a complex number, respectively. It can be shown that if $M \subseteq \mathbb{K}$ is a free $\mathbb{Z}$-module of rank $n$ with $\mathbb{Z}$-basis $\{x_1, \ldots,x_n\}$, then $\sigma_{\mathbb{K}}(M)$ is a full-rank lattice in $\mathbb{R}^n$ with basis $\{\sigma_{\mathbb{K}}(x_1), \ldots,\sigma_{\mathbb{K}}(x_n)\}$. Additionally, a  Gram matrix for \( \sigma_{\mathbb{K}}(M) \) is given by
\begin{equation} \label{Gram}
    G = \frac{1}{2} \left[ \Tr_{\mathbb{K}/\mathbb{Q}} (x_i \overline{x_j})\right]_{1 \leq i, j \leq n},
\end{equation}
where $\Tr_{\mathbb K/\mathbb{Q}}$ denotes the field trace, i.e., the sum of the 
$n$ embeddings of $\mathbb K.$ 

We denote by $\anelok$ the ring of integers of $\mathbb K$, which is the set of all algebraic integers in $\mathbb K$. If $M$ is a $\mathbb Z$-submodule of $\anelok$ of rank $n$, then the minimum Euclidean distance and volume of $\sigma_{\mathbb{K}}(M)$ are given, respectively, by
\[
2^{-\frac{1}{2}}\min_{\substack{x \in M \setminus \{0\}}} \sqrt{\Tr_{\mathbb K/\mathbb{Q}}(x\overline{x})} \quad \text{and} \quad 2^{-\frac{n}{2}}\sqrt{|\Disc(\mathbb K)|} \cdot [\mathcal{O}_{\mathbb K} : M],
\]
where $\Disc(\mathbb K)$ represents the discriminant of the number field $\mathbb K$. Therefore, since $\rho$ is half the minimum distance in $\Lambda$, the center density of $\sigma_{\mathbb K}(M)$ is given by
\begin{equation}\label{eq}
    \delta(\sigma_{\mathbb K}(M)) = \frac{\left(\min_{\substack{x \in M \\ x \ne 0}} \Tr_{\mathbb K/\mathbb{Q}}(x\overline{x})\right)^{\frac{n}{2}}}{2^n\sqrt{|\Disc(\mathbb K}| \cdot [\mathcal{O}_K : M]}.
\end{equation}

\section{Lattice $D_4$ via $\mathbb Q\left(\sqrt{ 2},\sqrt{- q}\right)$}
Let $\mathbb K = \mathbb Q\left(\sqrt{ 2},\sqrt{- q}\right)$ be a totally imaginary biquadratic field, where $q$ is a positive  square-free integer congruent to $3$ modulo $8$. Since $(-q) \equiv 1 \pmod 4,$ the set $\{\alpha_1,\alpha_2,\alpha_3,\alpha_4\}$ is an integral basis for $\mathbb{K}$, with
\[
\alpha_1 = 1, \quad \alpha_2 = \sqrt{2}, \quad \alpha_3 = \frac{1+i\sqrt{q}}{2}, \quad \alpha_4 = \frac{\sqrt{2}+i\sqrt{2q}}{2}.
\]
and
\[
 \Disc(\mathbb K) = 64q^2
\]
(see \cite{Williams}). Given $x = \sum_{i=1}^4 x_i \alpha_i \in \anelok$, where $x_1,x_2,x_3,x_4 \in \mathbb Z,$ we have:
\begin{align*}
    x &= x_1 + x_2\sqrt{2} + \frac{x_3}{2} + \frac{x_4\sqrt{2}}{2} + i\sqrt{q}\left(\frac{x_3 + x_4\sqrt{2}}{2}\right).
\end{align*}
Thus,
\begin{align*}
    x\overline{x} &= \left(x_1 + x_2\sqrt{2} + \frac{x_3}{2} + \frac{x_4\sqrt{2}}{2}\right)^2 + q\left(\frac{x_3 + x_4\sqrt{2}}{2}\right)^2 \\
    &= x_1^2 + 2x_2^2 + \frac{x_3^2}{4} + \frac{x_4^2}{2}  + x_1x_3 
     + 2x_2x_4 + q\left(\frac{x_3^2 + 2x_4^2}{4}\right)\\&\quad+ \sqrt{2}\left(2x_1x_2+x_1x_4+ x_2x_3+ \frac{(1+q)x_3x_4}{2}\right).
\end{align*}
We conclude that:
\begin{equation}\label{tracexx}
    \Tr_{\mathbb{K}/\mathbb{Q}}(x\overline{x}) = (2x_1 + x_3)^2 + 2(2x_2 + x_4)^2 + q(x_3^2 + 2x_4^2).
\end{equation}

\begin{prop}{}{}\label{caracterizacao-do-modulo} Let $j$ be a integer such that:
\[
\mathcal{M}_{j,q} := \left\{x = \sum_{i=1}^4 x_i \alpha_i \in \mathcal{O}_{\mathbb{K}} \mid x_3 \equiv -2x_1 + j(2x_2 + x_4) \mod{q}\right\}.
\]
Then,  $\mathcal{M}_{j,q}$ has rank $4$ and index $q$. \end{prop} 

\begin{proof} Since $x_3 \equiv -2x_1 + j(2x_2 + x_4) \pmod{q}$, there exists $t \in \mathbb Z$ such that $x_3 = -2x_1 + j(2x_2 + x_4) + tq$. Therefore, 
\begin{align*}
    x &= x_1\alpha_1 + x_2\alpha_2 + (-2x_1 + j(2x_2 + x_4) + tq)\alpha_3 + x_4\alpha_4 \\
    &= x_1(\alpha_1 - 2\alpha_3) + x_2(\alpha_2 + 2j\alpha_3) + t(q\alpha_3) + x_4(j\alpha_3 + \alpha_4).
\end{align*} It is straightforward to check that \( \{\alpha_1 - 2\alpha_3, \alpha_2 + 2j\alpha_3, q\alpha_3, j\alpha_3 + \alpha_4\} \) is linearly independent over \( \mathbb{Z} \). Moreover:
$$[\mathcal{O}_{\mathbb{K}} : \mathcal{M}_{j,q}] = 
\begin{vmatrix}
1 & \quad 0 & \quad -2 & \quad 0 \\
0 & \quad 1 & \quad 2j & \quad 0 \\
0 & \quad 0 & \quad q & \quad 0 \\
0 & \quad 0 & \quad j & \quad 1
\end{vmatrix} = q,$$
see \cite[Ch. 1, Theorem 1.17]{Stewart}. 
 \end{proof}

Additionally, consider the following auxiliary proposition:

\begin{prop}
Let \( q = p_1 \cdots p_m \), where \( m \geq 1 \) and \( p_j \) are distinct primes congruent to either \( 1 \) or \( 3 \) modulo \( 8 \), for \( j = 1, \dots, m \). Then, \(-2\) is a quadratic residue modulo \( q \).
\end{prop}
\begin{proof} By \cite[Ch. VI, Theorems 82 and 95]{Hardy}, for all primes \( p \equiv 1 \pmod{8} \) and \( p \equiv 3 \pmod{8} \), we have $-2$ is a quadratic residue modulo $p$. Indeed,
\[ \left(\frac{-2}{p}\right) = \left(\frac{-1}{p}\right)\left(\frac{2}{p}\right) = 1.\]
It then follows from \cite[pp. 287–289]{Stewart} that \( -2 \) is a quadratic residue modulo $ q = p_1 \cdots p_m$.

\end{proof}

This allows us to establish the following result:

\begin{prop}\label{lema1}
Let $q \equiv 3 \pmod 8$ be a positive square-free integer with all primes factors congruent to either $1$ or $3$ modulo $8,$ and let $j$ be a integer such that $j^2 \equiv -2 \pmod q$. If $x$ is a nonzero element of $$\mathcal{M}_{j,q} = \left\{x =\sum_{i=1}^4 x_i \alpha_i \in \mathcal{O}_{\mathbb{K}} \mid x_3 \equiv -2x_1 + j(2x_2 + x_4) \mod{q}\right\},$$ then $\Tr_{\mathbb{K}/\mathbb{Q}}(x\overline{x}) \geq 4q$.
\end{prop}

\begin{proof} 
First, note that since \( x_3 \equiv -2x_1 + j(2x_2 + x_4) \pmod{q}, \)
\[
(2x_1 + x_3)^2 \equiv -2(2x_2 + x_4)^2 \pmod{q}.
\]
Thus, there exists an integer \( k \geq 0 \) such that
\[
(2x_1 + x_3)^2 = -2(2x_2 + x_4)^2 + kq,
\]
and then, by \eqref{tracexx},
\begin{align*}
    \Tr_{\mathbb{K}/\mathbb{Q}}(x\overline{x}) = q(k + x_3^2 + 2x_4^2).
\end{align*}
Note that if \( |x_3| \geq 2 \) or \( |x_4| \geq 2 \), the result follows immediately. We now consider the following cases:

\begin{enumerate}
    \item[1.] If $x_3 = x_4 = 0$, it follows that
    \[
    \Tr_{\mathbb{K}/\mathbb{Q}}(x\overline{x}) = qk.
    \]
    If \( k = 0 \), then \( 4x_1^2 = -8x_2^2 \), which is only possible if \( x = 0 \), a contradiction. If \( k = 1 \), then
    \[
    q = 4x_1^2 + 8x_2^2 \equiv 0 \pmod{4},
    \]
    which contradicts \( q \equiv 3 \pmod{4} \). If \( k = 2 \), since \( \gcd(2, q) = 1 \),
    \[
    q = 2x_1^2 + 4x_2^2 \equiv 0 \pmod{2},
    \]
    again a contradiction. If \( k = 3 \), then
    \[
    3q = 4x_1^2 + 8x_2^2 \equiv 0 \pmod{4},
    \]
    equating an odd number to an even one, which is impossible. Therefore, \( \Tr(x\overline{x}) \geq 4q \) in this case.

    \item[2.] If $x_3 = 0$ and $ x_4 = \pm 1$, then
    \[
    \Tr_{\mathbb{K}/\mathbb{Q}}(x\overline{x}) = (k + 2)q.
    \]
    If \( k = 0 \), then \( 4x_1^2 = -2(2x_2 \pm 1)^2 \), which is impossible. If \( k = 1 \), then
    \[
    q = 8x_2^2 \pm 8x_2 + 2 + 4x_1^2 \equiv 2 \pmod{4},
    \]
    which contradicts \( q \equiv 3 \pmod{4} \). Thus, \( k \geq 2 \) and, consequently, \( \Tr_{\mathbb{K}/\mathbb{Q}}(x\overline{x}) \geq 4q \).

    \item[3.] If $x_3 = \pm 1$ and $ x_4 = 0$, it follows that
    \[
    \Tr_{\mathbb{K}/\mathbb{Q}}(x\overline{x}) = (k + 1)q.
    \]
    If \( k = 0 \), then \( (2x_1 \pm 1)^2 = -8x_2^2 \), a contradiction. If \( k = 1 \) or \( 2 \), then \( q \) or \(2q \equiv 1 \pmod{4} \), again a contradiction. Therefore, \( k \geq 3 \) and \( \Tr_{\mathbb{K}/\mathbb{Q}}(x\overline{x}) \geq 4q \).

    \item[4.] Finally, if $x_3 = \pm 1$ and $ x_4 = \pm 1$, we obtain  
    \[
    \Tr_{\mathbb{K}/\mathbb{Q}}(x\overline{x}) = (k + 3)q.
    \]
    Since \( k = 0 \) leads to a contradiction, we have \( k \geq 1 \), and thus  \newline \( \Tr_{\mathbb{K}/\mathbb{Q}}(x\overline{x}) \geq 4q \). 
\end{enumerate}

\end{proof} 

Finally, we can establish the main result of this section:

\begin{thm}\label{Teoremaprincipal}
Let $q \equiv 3 \pmod 8$ be a positive square-free integer with all primes factors congruent to either $1$ or $3$ modulo $8,$ and let $j$ be a integer such that $j^2 \equiv -2 \pmod q$. Thus, $\sigma_{\mathbb K}(\mathcal{M}_{j,q})$ is equivalent to the $D_4$-lattice.
\end{thm}

\begin{proof}
Since \([ \mathcal{O}_{\mathbb{K}} : \mathcal{M}_{j,q} ] = q\) and $\Disc(\mathbb K) = 64q^2,$ using \eqref{eq} and Proposition \ref{lema1}, we have:
\[
\delta\left(\sigma_{\mathbb{K}}\left(\mathcal{M}_{j,q}\right)\right) = \frac{\left(\min_{\substack{x \in M \\ x \ne 0}} \Tr_{\mathbb{K}/\mathbb{Q}}(x\overline{x})\right)^2}{2^4 \cdot \sqrt{\Disc(\mathbb{K})} \cdot [\mathcal{O}_{\mathbb{K}} : \mathcal{M}_{j,q}]} \geq \frac{(4q)^2}{16 \cdot 8q \cdot q} = \frac{1}{8}.
\]
As shown in \cite{Korkine}, the lattice $D_4$ has center density \( \frac{1}{8} \), which is the maximum possible among lattices in $\mathbb R^4$. Thus, $$\delta\left(\sigma_{\mathbb{K}}\left(\mathcal{M}_{j,q}\right)\right) = \frac{1}{8}.$$ Moreover, since the densest 4-dimensional lattice packing is unique (see \cite{Conway}, p. 15), the lattices of the form \( \sigma_{\mathbb{K}}(\mathcal{M}_{j,q}) \) are equivalent to \( D_4 \).

\end{proof}

By Dirichlet's Theorem \cite[Ch. II, Theorem 15]{Hardy}, there are infinitely many primes congruent to either $1$ or $3$ modulo $8$. Thus, there are infinitely many integers $q$ that satisfy the conditions of Theorem \ref{Teoremaprincipal}.

\section{Lattices $D_8$ and $E_8$ via $\mathbb Q\left(\sqrt{ 2},\sqrt{q},i\right)$}
In this section, consider $\mathbb L = \mathbb Q\left(\sqrt{ 2},\sqrt{ q},i\right)$ a totally imaginary triquadratic field, where $q$ has the same properties as in Theorem \ref{Teoremaprincipal}.

\begin{prop}\label{prod4d4}
    Let $\{z_1, z_2, z_3, z_4\}$ be a $\mathbb{Z}$-basis for $\mathcal{M}_{j,q}$  and $\mathcal{S}$ be the $\mathbb{Z}$-module generated by $\{z_1, z_2, z_3, z_4, iz_1, iz_2, iz_3, iz_4\}.$ Then, $\sigma_{\mathbb L}(\mathcal{S})$ is equivalent to $D_4 \oplus D_4.$
\end{prop}
\begin{proof}
     Since $z_i \in \mathbb{K} = \mathbb Q\left(\sqrt{ 2},\sqrt{ - q}\right)$ for all $1\leq j,k \leq 4,$  using trace properties: 
\begin{align*}
    &\Tr_{\mathbb{L}/\mathbb{Q}}(z_j\overline{z_k}) = 2 \Tr_{\mathbb{K}/\mathbb{Q}}(z_j\overline{z_k}),\\
    &\Tr_{\mathbb{L}/\mathbb{Q}}(iz_j\overline{iz_k}) = \Tr_{\mathbb{L}/\mathbb{Q}}(z_j\overline{z_k}) = 2 \Tr_{\mathbb{K}/\mathbb{Q}}(z_j\overline{z_k}) \, \mbox{ and } \\
    &\Tr_{\mathbb{L}/\mathbb{Q}}(iz_j\overline{z_k}) = \Tr_{\mathbb{L}/\mathbb{Q}}(z_j\overline{iz_k}) = 0.
\end{align*}
Then, by \eqref{Gram}, the Gram matrix of $\sigma_{\mathbb L}(\mathcal{S})$ is given by
\begin{equation*}
    G = 2\begin{bmatrix}
    G' &  0 \\
    0 & G'
\end{bmatrix},
\end{equation*}
where $G'$ is a Gram matrix of $\sigma_{\mathbb{K}}(\mathcal{M}_{j,q})$ derived from  $\{z_1,z_2,z_3,z_4\}$. Thus, $\sigma_{\mathbb L}(\mathcal{S})$ is equivalent to $D_4 \oplus D_4.$
\end{proof}
Classical Gram matrices, respectively, for the lattices  $D_4$ and $D_4 \oplus D_4$ are
\begin{equation} \label{GramD4}
\begin{aligned}
    G_1 = \begin{bmatrix}
2 & 0 & -1 & 0 \\
0 & 2 & -1 & 0 \\
-1 & -1 & 2 & -1\\
0 & 0 & -1 & 2 
\end{bmatrix} \quad \text{and} \quad G_2 = \begin{bmatrix}
    G_1 &  0 \\
    0 &  G_1
\end{bmatrix}.
\end{aligned}   
\end{equation}

 Since \( \sigma_{\mathbb{K}}(\mathcal{M}_{j,q}) \) is equivalent to \( D_4 \), there exist \( \gamma \in \mathbb{R} \setminus \{0\} \) and a unimodular matrix \( U \) such that $G_1 = \gamma U G' U^t$ (see \cite{Conway}). Moreover, if
\begin{align} \label{conjuntow}
    \begin{bmatrix}
    w_1 \\ w_2 \\ w_3 \\ w_4
    \end{bmatrix}
    = U \cdot
    \begin{bmatrix}
    z_1 \\ z_2 \\ z_3 \\ z_4
    \end{bmatrix},
\end{align}
then \( G_1 \) and $G_2$ are Gram matrices associated with the $\mathbb{Z}$-bases \( \{w_1,w_2,w_3,w_4\}\) of $\mathcal{M}_{j,q}$ and \( \{w_1,w_2,w_3,w_4, iw_1, iw_2, iw_3, iw_4\}\) of $\mathcal{S}$, respectively, up to scale factors.

In a specific case, the following proposition provides an explicit description of \( U \) and, consequently, of \( \{w_1, w_2, w_3, w_4\} \).

\begin{prop}
    Let $q = 4k^2+4k+3,$ $j = 2k+1$ and 
    $$U = \begin{bmatrix}
        1 & \quad 0 & \quad 0 & \quad 0\\
        j & \quad 1 & \quad 0 & \quad 0\\
        -k & \quad -k - 1& \quad 1 & \quad 0\\
        -1 & \quad k & \quad -1 & \quad 1
    \end{bmatrix},$$
    where $k \in \mathbb{Z}$. 
    If we have $z_1= \alpha_1 - 2\alpha_3$, $z_2 = \alpha_2 + 2j\alpha_3$, $z_3 = q\alpha_3$ and  $z_4 = j\alpha_3 + \alpha_4$,
    then the Gram matrix $G_{j,q}$ associated with $B = \{z_1,z_2,z_3,z_4\}$ satisfies 
    $G_1 =\frac{1}{q}UG_{j,q}U^t,$ where $G_1$ is given by \eqref{GramD4}.
\end{prop}
\begin{proof}
    First, note that $q \equiv 3 \pmod 8.$ Indeed, since we have that the number $l = k(k+1)/2$ is an integer, it follows that
    $$q = 4k^2+4k+3 = 4k(k+1)+3 = 8l+3 \equiv 3 \pmod 8.$$
    Moreover, we have that
    $$j^2 = (2k+1)^2 = 4k^2+4k+1 \equiv -2 \pmod q.$$
    Next, note that, by the proof of Proposition \ref{caracterizacao-do-modulo}, $B$ is a $\mathbb Z$-basis of $M_{j,q}$. Additionally, by direct computation, we obtain: 
    \begin{equation*}\begin{aligned}
        z_1 &= -i\sqrt{q}, \quad z_2 = j+\sqrt{2}+ji\sqrt{q},\\ \displaystyle z_3 &= \frac{q+qi\sqrt{q}}{2} \quad \text{and} \quad z_4 = \frac{j+\sqrt{2}+ji\sqrt{q}+i\sqrt{2q}}{2}.
    \end{aligned}     
    \end{equation*}
    Then, the Gram matrix of $\sigma_{\mathbb{K}}(\mathcal{M}_{j,q})$ associated to $B$ is
    $$G_{j,q} = q\begin{bmatrix}
        2 & \quad -2j & \quad -q& \quad -j\\
-2j& \quad \frac{2j^2q + 2j^2 + 4}{q}& \quad j(q + 1)& \quad \frac{j^2q + j^2 + 2}{q}\\
-q& \quad j(q + 1)& \quad \frac{q(q + 1)}{2}& \quad \frac{j(q + 1)}{2}\\
-j& \quad \frac{j^2q + j^2 + 2}{q}& \quad \frac{j(q + 1)}{2} &\quad \frac{(j^2+2)(q+1)}{2q} 
    \end{bmatrix}.$$
    Using that $q = 4k^2+4k+3$ and $j = 2k+1$, we obtain that $G_1 = \frac{1}{q}UG_{j,q}U^t$.
    \end{proof}

In the following theorems, consider that $G_1$ from \eqref{GramD4} is the Gram matrix derived from the set  \( \{w_1, w_2, w_3, w_4\} \) from \eqref{conjuntow}, up to scaling.

\begin{thm}\label{GeraD8}
    Let $v = \frac{i-1}{2}w_1+\frac{-i-1}{2}w_2-w_3-w_4.$ Then the set
    $$
    B_1 = \left\{w_1, w_2, w_3, w_4, v, iw_2, iw_3, iw_4\right\}
    $$
    is a basis for a $\mathbb{Z}$-module $\mathcal{M}$ such that $\sigma_{\mathbb{L}}(\mathcal{M})$ is a lattice equivalent to $D_8.$
\end{thm}

\begin{proof}
    From \eqref{Gram}, we have that, for all $1 \leq j,k \leq 4,$ 
    \begin{equation}\label{tracesd4}
        \begin{aligned}
            \Tr_{\mathbb{L}/\mathbb{Q}}(iw_j\overline{iw_k}) &= \Tr_{\mathbb{L}/\mathbb{Q}}(w_j\overline{w_k}) = 2g_{jk}, \\
        \Tr_{\mathbb{L}/\mathbb{Q}}(iw_j\overline{w_k}) &= 0,
        \end{aligned}
    \end{equation}
    where $g_{jk}$ are the entries of $G_1$ from \eqref{GramD4}, up to scaling. Moreover, given $u \in B_1$, it follows from the linearity properties of the trace that:
    \begin{align*}
        \Tr_{\mathbb{L}/\mathbb{Q}}(v\overline{u})
        &= \frac{1}{2} \Tr_{\mathbb{L}/\mathbb{Q}}(iw_1\overline{u}) - \frac{1}{2} \Tr_{\mathbb{L}/\mathbb{Q}}(w_1\overline{u}) 
        - \frac{1}{2} \Tr_{\mathbb{L}/\mathbb{Q}}(iw_2\overline{u}) \\
        &\quad - \frac{1}{2} \Tr_{\mathbb{L}/\mathbb{Q}}(w_2\overline{u}) - \Tr_{\mathbb{L}/\mathbb{Q}}(w_3\overline{u}) - \Tr_{\mathbb{L}/\mathbb{Q}}(w_4\overline{u}).
    \end{align*}
    Thus, using \eqref{tracesd4}, we compute:
    \[
    \begin{array}{|c|c|c|c|c|c|c|c|c|}
    \hline
   u  & w_1 & w_2 & w_3 & w_4 & v & iw_2 & iw_3 & iw_4 \\
    \hline
    \Tr_{\mathbb{L}/\mathbb{Q}}(v\overline{u}) & 0 & 0 & 0 & -2 & 4 & -2 & 0 & 0 \\
    \hline
    \end{array}
    \]
    Therefore, by \eqref{Gram}, since \(\Tr_{\mathbb{L}/\mathbb{Q}}(u\overline{v}) = \Tr_{\mathbb{L}/\mathbb{Q}}(v\overline{u})\), the Gram matrix, up to scaling, associated with \(B_1\) is:
    \begin{equation}\label{GramD8}
        \begin{bmatrix}
        2 & 0 & -1 & 0 & 0 & 0 & 0 & 0 \\
        0 & 2 & -1 & 0 & 0 & 0 & 0 & 0 \\
        -1 & -1 & 2 & -1 & 0 & 0 & 0 & 0 \\
        0 & 0 & -1 & 2 & -1 & 0 & 0 & 0 \\
        0 & 0 & 0 & -1 & 2 & -1 & 0 & 0 \\
        0 & 0 & 0 & 0 & -1 & 2 & -1 & 0 \\
        0 & 0 & 0 & 0 & 0 & -1 & 2 & -1 \\
        0 & 0 & 0 & 0 & 0 & 0 & -1 & 2 \\
        \end{bmatrix},
    \end{equation}
    which is a classical Gram matrix of \(D_8\).
\end{proof}

\begin{thm}\label{GeraE8}
    Let \( v_1 = w_1 - w_2, \quad v_2 = \frac{i-1}{2}w_1 + \frac{-i-1}{2}w_2 - w_3 - w_4, \) and  
    \[
    v_3 = (i+1)w_1 + (i+1)w_3 + \frac{i+1}{2}(w_2 + w_4).
    \]
    Then the set
    \[
    B_2 = \left\{v_1, w_2, w_3, w_4, v_2, iw_2, iw_3, v_3\right\}
    \]
    is a basis for a $\mathbb{Z}$-module $\mathcal{M}$ such that $\sigma_{\mathbb{L}}(\mathcal{M})$ is a lattice equivalent to $E_8.$
\end{thm}

\begin{proof}
    First, note that \( v_2 \) is the same element as \( v \) from Theorem \ref{GeraD8}.  
    Therefore, we only need to compute the entries of the Gram matrix involving \( v_1 \) and \( v_3 \).  
    For any \( u \in B_2 \), by the linearity properties of the trace, we have:
    \begin{align*}
        \Tr_{\mathbb{L}/\mathbb{Q}}(v_1\overline{u}) 
        &= \Tr_{\mathbb{L}/\mathbb{Q}}(w_1\overline{u}) - \Tr_{\mathbb{L}/\mathbb{Q}}(w_2\overline{u}), \\
        \Tr_{\mathbb{L}/\mathbb{Q}}(v_3\overline{u}) 
        &= \Tr_{\mathbb{L}/\mathbb{Q}}(iw_1\overline{u}) + \Tr_{\mathbb{L}/\mathbb{Q}}(w_1\overline{u}) + \frac{1}{2} \Tr_{\mathbb{L}/\mathbb{Q}}(iw_2\overline{u}) \\
        &\quad + \frac{1}{2} \Tr_{\mathbb{L}/\mathbb{Q}}(w_2\overline{u}) + \Tr_{\mathbb{L}/\mathbb{Q}}(iw_3\overline{u}) + \Tr_{\mathbb{L}/\mathbb{Q}}(w_3\overline{u}) \\
        &\quad + \frac{1}{2} \Tr_{\mathbb{L}/\mathbb{Q}}(iw_4\overline{u}) + \frac{1}{2} \Tr_{\mathbb{L}/\mathbb{Q}}(w_4\overline{u}).
    \end{align*}
    Thus, using \eqref{tracesd4}, we compute:
    \[
    \begin{array}{|c|c|c|c|c|c|c|c|c|}
    \hline
    u  & v_1 & w_2 & w_3 & w_4 & v_2 & iw_2 & iw_3 & v_3 \\
    \hline
    \Tr_{\mathbb{L}/\mathbb{Q}}(v_1\overline{u}) & 8 & -4 & 0 & 0 & 0 & 0 & 0 & 2 \\
    \hline
    \Tr_{\mathbb{L}/\mathbb{Q}}(v_3\overline{u}) & 2 & 0 & 0 & 0 & 0 & 0 & 0 & 4 \\
    \hline
    \end{array}
    \]
    Therefore, by \eqref{Gram}, since \(\Tr_{\mathbb{L}/\mathbb{Q}}(u\overline{v}) = \Tr_{\mathbb{L}/\mathbb{Q}}(v\overline{u})\), the Gram matrix, up to scaling, associated with $B_2$ is:
    \begin{equation}\label{GramE8}
        \begin{bmatrix}
            4 & -2 & 0 & 0 & 0 & 0 & 0 & 1 \\
            -2 & 2 & -1 & 0 & 0 & 0 & 0 & 0 \\
            0 & -1 & 2 & -1 & 0 & 0 & 0 & 0 \\
            0 & 0 & -1 & 2 & -1 & 0 & 0 & 0 \\
            0 & 0 & 0 & -1 & 2 & -1 & 0 & 0 \\
            0 & 0 & 0 & 0 & -1 & 2 & -1 & 0 \\
            0 & 0 & 0 & 0 & 0 & -1 & 2 & 0 \\
            1 & 0 & 0 & 0 & 0 & 0 & 0 & 2
        \end{bmatrix},
    \end{equation}
    which is a classical Gram matrix of \(E_8\).
\end{proof}



\begin{thebibliography}{99}

\bibitem{Andrade} A. A. Andrade, A. J. Ferrari, C. W. O. Benedito and S. I. R. Costa,
Constructions of algebraic lattices,
{\it Comput. Appl. Math.}, {\bf 29} (2010), 1-13.

\bibitem{Conway} J. H. Conway and N. J. A. Sloane,
{\it Sphere Packings, Lattices and Groups},
3rd Edition, Springer Verlag, New York (1999).

\bibitem{Trajano} A. L. Flores, J. C. Interlando, T. P. da Nóbrega Neto and A. L. Contiero,
A new number field construction of the lattice $E_8$,
{\it Beitr. Algebra Geom.}, {\bf 54} (2013), 503-508.

\bibitem{Eva2} E. Bayer-Fluckiger,
Ideal lattices,
In: {\it Proceedings of the Conference Number Theory and Diophantine Geometry}, Zurich, 1999,
Cambridge Univ. Press (2002), 168-184.

\bibitem{Eva} E. Bayer-Fluckiger,
Lattices and number fields,
{\it Contemp. Math.}, {\bf 241} (1999), 69-84.

\bibitem{Eva-Ivan} E. Bayer-Fluckiger and I. Suarez,
Ideal lattices over totally real number fields and Euclidean minima,
{\it Arch. Math.}, {\bf 86}(3) (2006), 217-225.

\bibitem{Boutros} J. Boutros, E. Viterbo, C. Rastello and J. C. Belfiori,
Good lattice constellations for both Rayleigh fading and Gaussian channels,
{\it IEEE Trans. Inform. Theory}, {\bf 42}(2) (1996), 502-517.

\bibitem{Craig} M. Craig,
A cyclotomic construction of the Leech's lattice,
{\it Mathematika}, {\bf 25} (1978), 236-241.

\bibitem{Hardy} G. H. Hardy and E. M. Wright,
{\it An Introduction to the Theory of Numbers},
6th Edition, Oxford University Press, Oxford (2008).

\bibitem{Interlando} J. C. Interlando, J. O. D. Lopes and T. P. da Nóbrega Neto,
A new number field construction of the $D_4$ lattice,
{\it Int. J. Appl. Math.}, {\bf 31}(2) (2018), 299-305.

\bibitem{Jorge} G. C. Jorge, A. J. Ferrari and S. I. R. Costa,
Rotated $D_n$ lattices,
{\it J. Number Theory}, {\bf 132} (2012), 2397-2406.

\bibitem{Jorge2} G. C. Jorge and S. I. R. Costa,
On rotated $D_n$-lattices constructed via totally real number fields,
{\it Arch. Math.}, {\bf 100} (2013), 323-332.

\bibitem{Jorge-Densos} G. C. Jorge, A. A. de Andrade, S. I. R. Costa and J. E. Strapasson,
Algebraic constructions of densest lattices,
{\it J. Algebra}, {\bf 429} (2015), 218-235.

\bibitem{Korkine} A. Korkine and G. Zolotareff,
Sur les formes quadratique positive quaternaires,
{\it Math. Ann.}, {\bf 5} (1872), 581-583.

\bibitem{Samuel} P. Samuel,
{\it Algebraic Theory of Numbers},
Hermann, Paris (1970).

\bibitem{Shannon} C. E. Shannon,
A mathematical theory of communication,
{\it Bell Syst. Tech. J.}, {\bf 27} (1948), 379-423.

\bibitem{Stewart} I. Stewart,
{\it Algebraic Number Theory and Fermat's Last Theorem},
4th Edition, Chapman and Hall/CRC, New York (2015).

\bibitem{Viazovska} M. S. Viazovska,
The sphere packing problem in dimension 8,
{\it Ann. Math.}, {\bf 185}(3) (2017), 991-1015.

\bibitem{Williams} K. S. Williams,
Integers of biquadratic fields,
{\it Can. Math. Bull.}, {\bf 13}(4) (1970), 519-526.

\end{thebibliography}
\end{document}